\documentclass[11pt]{article}
\usepackage{amsfonts}
\usepackage{amssymb,amsmath,amsthm,amscd}
\usepackage{epsfig}
\usepackage{graphicx,hyperref,url}
\usepackage[notref,notcite,final]{showkeys}

\newtheorem{theo}{Theorem}
\newtheorem{prop}{Proposition}
\newtheorem{ques}{Question}
\newtheorem{defi}{Definition}
\newtheorem{lema}{Lemma}
\newtheorem{rema}{Remark}

\newcommand{\Cset}{\mathbb{C}}

\newcommand{\Rset}{\mathbb{R}}

\newcommand{\fm}{\phantom{-}}
\newcommand{\df}{\displaystyle\frac}
\newcommand{\dd}{\mbox{ d}}
\newcommand{\fpn}{\frac{\pi}{n}}
\newcommand{\eps}{\varepsilon}

\begin{document}

\Large{

\vspace*{-5mm}

\begin{center}

{\LARGE Really perverse periodic solutions of the planar $N$-body problem}
\end{center}

\begin{center}{Dominique Bang$^{(1)}$, Alain Chenciner$^{(2,3)}$, Carles Sim\'o$^{(4,5)}$}
\end{center}

\vspace*{1mm}

{\small

\noindent$^{(1)}$ Bank Of America,
51 Rue La Boétie, 75008 Paris, France,

\noindent $ $\hspace*{3mm} dominique.bang@bofa.com

 \noindent $^{(2)}$ Astronomie et Syst\`emes Dynamiques, LTE, 77, avenue Denfert-Rochereau,
75014 Paris, France

\noindent $ $\hspace*{3mm}  alain.chenciner@obspm.fr

\noindent $^{(3)}$ UFR de Math\'ematiques, Universit\'e Paris Cit\'e

\noindent $^{(4)}$ Departament de Matem\`atiques i Inform\`atica, Univ. de
Barcelona, Gran Via, 585, Barcelona 08007, Spain,

Honorary Professor Emeritus

\noindent $ $ \hspace*{3mm} carles@maia.ub.es

\noindent $^{(5)}$ Reial Acad\`emia de Ci\`encies i Arts de Barcelona,
La Rambla, 115, Barcelona 08002, Spain

\vspace*{1cm}

}

ABSTRACT. Examples are given of solutions of the planar $N$-body problem which remain the same for at least two systems of masses with the same sum and same center of mass. The least value of $N$ achieved up to now with this property is 474, a number which had been announced in the first author's thesis.

\noindent N.B. Because of renewed interest, we choose to make public after a slight revision this paper written in 2002 but never published.

\vspace*{6mm}
\begin{center} THE PROBLEM AND MOTIVATIONS \end{center}

\vspace*{1mm}

Let $X(t)=(\vec r_1(t),\vec r_2(t),\ldots,\vec r_N(t))$ be a planar
solution of the Newtonian $N$-body problem with masses $m_1,m_2,\ldots,m_N$.  
The following questions were posed in \cite{C}:

\begin{ques} \label{q1} Does there exist $(m'_1,m'_2,\ldots,m'_N)$, another
system of masses, for which $X(t)$ is still a solution? \end{ques}

\begin{ques} \label{q2} Same as Question \ref{q1} but insisting that the sum $M=
\sum_{i=1}^N{m_i}$ of the masses and the center of mass $\vec r_G=\df{1}{M} 
\sum_{i=1}^N{m_i\vec r_i}$ do not change. \end{ques}

\begin{defi} If the answer to the first (resp. second) question is yes, we shall
say $X(t)$ is a {\rm perverse} (resp. {\rm really perverse}) solution and the 
allowed systems of masses will be called {\rm admissible}. \end{defi}

\vspace*{2mm}

Beyond the intrinsic interest of these natural questions there is another 
motivation. Recall \cite{Secm,CGMS} that a planar {\em choreography} is a
periodic solution ${\cal C}(t)=\left(q(t+\frac{T}{N}),\ldots, q(t+\frac{(N-1)T}
{N}),\,q(t+T)=q(t)\right)$ of the $N$-body problem such that all $N$ bodies
follow the same closed plane curve $q(t)$ with equal time spacing.

The question of the existence of choreographies whose masses are not all equal
is at the origin of the notion of perverse solution. Indeed, it is noticed in 
\cite{C2} that if such a choreography exists, it is a perverse choreography: 
masses all equal to the mean $\frac{M}{N}$ form an admissible set of masses.
Moreover, any perverse choreography is really perverse.
Partial results have been given in \cite{C}. In particular one has

\begin{prop} The perverse solutions of the planar 3-body problem are exactly 
the collinear homographic solutions.\end{prop}

\begin{prop} If $N\leq 4$, the planar $N$-body problem does not possess any 
really perverse solution. \end{prop}

\vspace*{2mm}

Our present goal is to show that {\em really perverse solutions exist}.
The simplest periodic solutions of the planar $N$-body problem being given by {\em relative equilibria}, we look for {\em really perverse relative equilibria}.
\smallskip

A simple configuration can be obtained by putting a central mass at the center
of masses and additional masses at the vertices of $k$ regular concentric 
$n$-gons, and then $N=1+kn.$ All the masses at the vertices of a given $n$-gon
are supposed to be equal (in case of a single regular $n$-gon this is a
necessary condition, see \cite{PW}), but different $n$-gons can have different
masses.

Other relative equilibria giving rise to really perverse solutions will also be
studied and generalizations to different potentials, other than Newtonian one,
will be discussed.

\vspace*{2mm}

\begin{center} STATEMENT OF RESULTS \end{center}

\vspace*{2mm}

\noindent We study several relative equilibria. In particular we prove the 
existence of {\em really perverse relative equilibria} in three cases:
\begin{itemize}
\item [1)] a central mass surrounded by 2 regular $n$-gons which differ
angularly by $\pi/n$, for $n\ge 237$,
\item [2)] a central mass surrounded by a regular $n$-gon and a regular
$2n$-gon, for $n\geq 237,$
\item [3)] a central mass surrounded by three homothetic regular $n$-gons,
for $456\leq n\leq 874$.
\end{itemize}

\vspace*{1mm}

\noindent The second case may be thought of as a central mass surrounded by 
three $n$-gons, two of them with the same radius and obtained one from the 
other by a rotation of angle $\pi/n$. See the next figure.
The proof is entirely analytical in the second case, partly numerical
in the first and the third ones.

\begin{center}
 \includegraphics[scale=1.05]{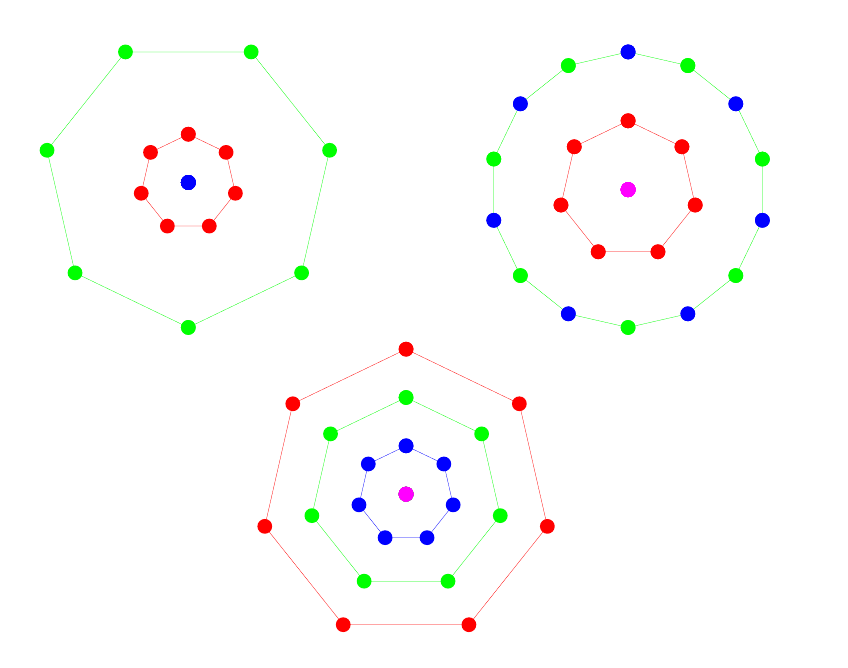}
\end{center}

\vspace*{2mm}

\begin{center} THE EQUATIONS \end{center}

\vspace*{1mm}

Let $m_0$ be the mass of the central body and $m_j$ the masses on the $j$-th
polygon.

\vspace*{1mm}

Identifying $\Rset^2$ with $\Cset$, let $x_l^j=\rho_j\exp(i(\frac{2\pi l}{n}+
\theta_j))$ be the location of the $l$-th vertex of the $j$-th polygon.
We shall assume $\rho_1\le \rho_2\le\ldots\le\rho_{k-1}\le\rho_k.$ As a 
normalization we can assume $\rho_k=1$ and $\theta_k=0.$ Furthermore the
values of $\theta_j,\; j=1,\ldots,k-1$ are assumed to be either 0 or
$\pi/n$ modulus $2\pi/n$, the only angles producing forces pointing to the
center of masses.

Note that in the plots shown in the previous figure, for aesthetic reasons these conditions are satisfied assuming that the angles are measured starting at zero on the positive part of the vertical axis.
\goodbreak

\vspace*{1mm}

It is not restrictive to assume that the angular velocity of the relative
equilibrium is $\omega=1.$
Then, for the Newtonian case, the equations of a relative equilibrium motion
$x_l^j(t)=x_l^je^{it}$ are

$$m_j\ddot x_m^j(t)=-\frac{m_jm_0}{\rho_j^3}x_m^j+
\sum_{s=1}^km_s\sum_{l=1}^{n^*(s/j)}m_j\frac{x_{m+l}^s-x_m^j}{|x_{m+l}^s-x_m^j|^3},$$

\noindent that is

$$\rho_j=m_0\frac{1}{\rho_j^2}+\sum_{s=1}^km_s{\sum_{l=1}^{n^*(s/j)}
\frac{\rho_j-\rho_s\cos(\frac{2\pi l}{n}+\theta_s-\theta_j)}
{\left(\rho_j^2+\rho_s^2-2\rho_j\rho_s\cos(\frac{2\pi l}{n}+\theta_s-\theta_j)
\right)^{3/2}}},$$
for $j=1,\ldots,k,$  where

\vspace*{-1mm}

$$n^*(z)=n\;{\mbox { if }}\; z\not=1,\, n^*(z)=n-1\;{\mbox { if }}\; z=1.$$

\noindent This can be written as
$$m_0+\sum_{s=1}^k{m_sH_n(\rho_s/\rho_j,\theta_s-\theta_j)}=\rho_j^{3},$$
where

\begin{equation*}
H_n(x,\theta) = \left\{
\begin{aligned}
&\sum_{l=1}^{n-1}\frac{1-\cos(\frac{2\pi l}{n})}{\left(1+1-2\cos(\frac{2\pi l}{n})\right)^{3/2}}=\frac{1}{4}\sum_{l=1}^{n-1}\frac{1}{\sin(\frac{\pi l}{n})}\quad \text {if $x=1, \theta=0$ (mod $\frac{2\pi}{n}$),}\\
&\sum_{l=1}^{n}\frac{1-x\cos(\frac{2\pi l}{n}+\theta)}{\left(1+x^2-2x\cos(\frac{2\pi l}{n}+\theta)\right)^{3/2}}\quad \text{otherwise}.
\end{aligned}
\right.
\end{equation*}
{\it Note that when $x=1$, that is $\rho_j=\rho_s$, the only values of $\theta$ at stake are $\theta=0$ if $j=s$ or $\theta =\frac{\pi}{n}$ mod $\frac{2\pi}{n}$ if $j\not=s$.}
 
\smallskip

\noindent If we insist that the total mass $M$ be fixed, the $k+1$ masses
$m_0,m_1,\ldots,m_k$ must satisfy a system of $k+1$ linear equations:

\vspace*{-1mm}

\begin{equation*}
\begin{pmatrix}1&H_n(1,0)&\cdots&H_n(\rho_k/\rho_1,\theta_k-\theta_1)\\
1&H_n(\rho_1/\rho_2,\theta_1-\theta_2)&\cdots &H_n(\rho_k/\rho_2,\theta_k-\theta_2)\\
\cdot&\cdots&\cdots&\cdots\\
1&H_n(\rho_1/\rho_k,\theta_1-\theta_k)&\cdots&H_n(1,0)\\
1&n&\cdots&n\end{pmatrix}
\begin{pmatrix}m_0\\ m_1\\ \cdot \\ \cdot \\ m_k \end{pmatrix}=\begin{pmatrix}\rho_1^3\\ \rho_2^3\\ \cdot\\
\rho_k^3\\ M
\end{pmatrix}
\end{equation*}    

It is convenient to introduce the notations

$$h_n(x,\theta)=\frac{1}{n}H_n(x,\theta),\quad
  k_n(x,\theta)=h_n(x,\theta)-1,\quad \delta_n=k_n(1,0),$$

and to write the equations as follows

$$\hspace*{-4mm}(*)\;\;\left(\!\! \begin{array}{ccc}
\delta_n&\cdots& k_n(\rho_k/\rho_1,\theta_k\!-\!\theta_1)\\
k_n(\rho_1/\rho_2,\theta_1\!-\!\theta_2)&\cdots &
k_n(\rho_k/\rho_2,\theta_k\!-\!\theta_2)\\ \cdots&\cdots&\cdots\\
k_n(\rho_1/\rho_k,\theta_1\!-\!\theta_k)&\cdots&\delta_n \end{array}\!\!\right)
\!\!\left(\!\!\begin{array}{c}nm_1\\ \cdots\\ \cdots\\nm_k\end{array}\!\right)
\!\!=\!\!
\left(\!\!\begin{array}{c}\rho_1^3\!-\!M\\ \rho_2^3\!-\!M\\ \cdots\\ 
\rho_k^3\!-\!M \end{array}\!\!\right)$$

\vspace*{2mm}

Hence, the problem is first to show the existence of $(k-1)$-uples
$\rho_1,\ldots,\rho_{k-1}$ (remember $\rho_k=1$) and $(k-1)$-uples 
$\theta_1,\ldots,\theta_{k-1}$ (remember $\theta_k=0$) such that the
determinant of the system vanishes, second to check that the affine subspace 
of solutions $m_0,m_1,\ldots,m_k$ is not empty and meets the positive quadrant.

\vspace*{2mm}

The above $k\times k$ matrix in $(*)$ will be denoted by $A_{n;k}=A$ .

\vspace*{2mm}

\begin{center} SOME PROPERTIES OF $\delta_n$ \end{center}

\vspace*{1mm}

We shall use the following asymptotic representation of $\delta_n$ (for derivation and additional information see \cite{MS} where
$\delta_n=\frac{A(n)}{n}-1$):

$$\begin{array}{rcl}\delta_n&=&
\displaystyle{\frac{1}{4n}{\sum_{l=1}^{n-1}\frac{1}{\sin(\frac{\pi l}{n})}}-1}\\
&=& \displaystyle{\frac{1}{2\pi}\left(\gamma+\log\frac{2n}{\pi}\right)+
{\sum_{k\ge 1}\frac{(-1)^k(2^{2k-1}-1)B^2_{2k}\pi^{2k-1}}{(2k)(2k)!}
\frac{1}{n^{2k}}}}-1, \end{array}$$

\noindent where $\gamma\!\approx\!0.577215664901532860606512090082$ stands for
the Euler-Mascheroni constant and $B_{2k}$ for the Bernoulli numbers. A few
terms give

\vspace*{1mm}

$$\delta_n\!\approx\! \displaystyle{\frac{1}{2\pi}\!\!\left(\!\gamma\!+\!\log\!
\frac{2n}{\pi}\!\right)\!\!-\!\frac{\pi}{144n^2}\!+\!\frac{7\pi^3}{86400n^4}\!-
\!  \frac{31\pi^5}{7620480n^6}\! +\!\frac{127\pi^7}{290304000n^8}}-\!1, $$

\vspace*{3mm}

\noindent which produce errors below $10^{-12}$ for $n\ge 9$ (and even below 
$1.5\times 10^{-6}$ for $n\ge 2$).

\vspace{4mm}

\noindent In particular, considering that $\delta_n$ is strictly increasing
(see \cite{BE}, Proposition 10), it follows from the asymptotic expansion
(which makes sense for $n\in\Rset$) that $\delta_n<0$ (resp. $\delta_n>0$) if
$n \leq 472$ (resp. $n \geq 473$) and $\delta_{n}$ vanishes for
$n\approx 472.270995.$

\vspace*{3mm}

\begin{center} STUDY OF THE CASES OF $1$ OR $2$ POLYGONS \end{center}

\vspace*{3mm}

\noindent {\bf The case of a single polygon ($k=1$)}

\vspace*{5mm}

\noindent The matrix $A$ reduces to $A_{n;1}=(\delta_n)$ and from the 
properties of $\delta_n$ follows

\begin{prop} \label{pk1} No really perverse solutions exist if $k=1.$\end{prop}

\vspace*{6mm}

\noindent {\bf The case of two polygons ($k=2$)}

\vspace*{5mm}

\noindent We shall consider two subcases

\vspace*{3mm}

\noindent {\bf 1) 2 homothetic polygons (unshifted case)}

\begin{prop} \label{pk2h} A relative equilibrium of a central mass surround\-ed 
by two homothetic regular $n$-gons is never really perverse. \end{prop}

\vspace*{3mm}

\noindent {\bf Proof.} Let $x=\rho_1/\rho_2<1.$ The determinant of $A_{n;2}(x)$
becomes 
$$\det A_{n;2}(x)=\delta_n^2-k_n(x)k_n(x^{-1}),$$
where we used the shorthand notation $k_n(x)=k_n(x,0)$. 

\vspace*{2mm}

\noindent It is sufficient to show that $k_n(x)$ is positive when $0<x<1$ and 
negative when $x>1$.

\vspace*{2mm}

\noindent As $k_n(0)=0$ and $\lim_{x\to+\infty}k_n(x)=-1$, this follows from

\vspace*{3mm}

\begin{lema} \label{l1} For each integer $n$, the function $x\mapsto k_n(x)$
is strictly increasing in each of the intervals $[0,1[$ and $]1,+\infty[$ : it
increases from $0$ to $\infty$ in the first case, from $-\infty$ to $-1$ in the
second one.
\end{lema}

\noindent To prove Lemma \ref{l1} it is convenient to use the representation
of the potential 

$$V_n(x,\theta)=\frac{1}{n}{\sum_{j=1}^n\frac{1}{(1+x^2-2x\cos(\frac{2\pi j}{n}
+\theta))^{1/2}}}$$

\vspace*{2mm}

\noindent for $0\le x<1$ as

$$ V_n(x,\theta)\!=\!\frac{1}{\pi}\!\int_0^1\!\!\frac{1}{((1\!-\!\tau)\tau)
^{1/2}}\frac{1}{(1\!-\!x^2\tau)^{1/2}}\frac{1-(x\tau)^{2n}}{1\!+\!(x\tau)^{2n}
\!-\!2(x\tau)^n \cos(n\theta)}\dd \tau $$ 

\vspace*{2mm}

\noindent which is obtained from a general theory of averaging analytic functions
on the vertices of a regular polygon (see \cite{BE,B}). 

\vspace*{3mm}

\noindent Then, from the relation 
$$h_n(x,\theta)=\frac{\partial}{\partial x}(xV_n(x,\theta))$$
and the fact that if $\theta=0$ the coefficients of all the powers of $x$ in
$ V_n(x,0)$ are positive, as follows from the integral representation, Lemma
\ref{l1} is obtained.

\vspace*{2mm}

A similar proof works in the case $x>1.$

\vspace*{1mm}

\noindent {\bf 2) 2 shifted polygons}

\begin{prop} \label{pk2s} For $n\ge 237$, there exists a really perverse 
relative equilibrium whose configuration is a central mass surrounded by 2 
shifted regular $n$-gons sharing the same symmetry axes. \end{prop}

The angle of rotation (or shift) with respect to homothetic polygons is,
of course, $\pi/n.$

\vspace*{1mm}

\begin{rema} The case $n=237$ gives the example of a really perverse solution
with $N=475$ total bodies. As will be seen, it is possible to take $m_0=0$ if
$n<473.$ In that case one reduces $N$ to $474$, the smallest number of bodies 
known at present in the Newtonian case for a really perverse solution.\end{rema}

\vspace*{2mm}

\noindent {\bf Proof.} The determinant which must be zero becomes now

$$f(x):=\det A_{n;2}(x,\fpn)=\delta_n^2-k_n(x,\fpn)k_n(x^{-1},
\fpn).$$

\vspace*{2mm}

\noindent One has $f(0)=\delta_n^2>0.$ On the other hand $f(1)=\delta_n^2-
(k_n(1,\fpn))^2,$ and using Lemma \ref{lcolec} (see later), we obtain 

$$f(1)=-2(\delta_{2n}-\delta_n)\times 2\delta_{2n},$$

\vspace*{2mm}

\noindent whose sign is opposite to the one of $\delta_{2n},$ because 
$\delta_n$ increases with $n.$

\vspace*{4mm}

\noindent This ensures the existence of at least one solution $x$ if $n\ge 237.$

\vspace*{4mm}

\noindent One checks numerically that this solution is unique. The behavior of
the solution for $n$ large, of the form $\;x=1-c\sqrt{\log(n)-\log(237)}/n^2\;$ 
for some $c>0,$ suggests that an asymptotic analysis seems feasible. Furthermore
for $n>1000$ one has checked that using $c=3.905$ the relative errors in $1-x$ ,
given by this expression, are less than $0.001$. The figure shown below has been
done with the values of $x$ computed numerically.

\begin{center}
  \includegraphics[scale=0.5]{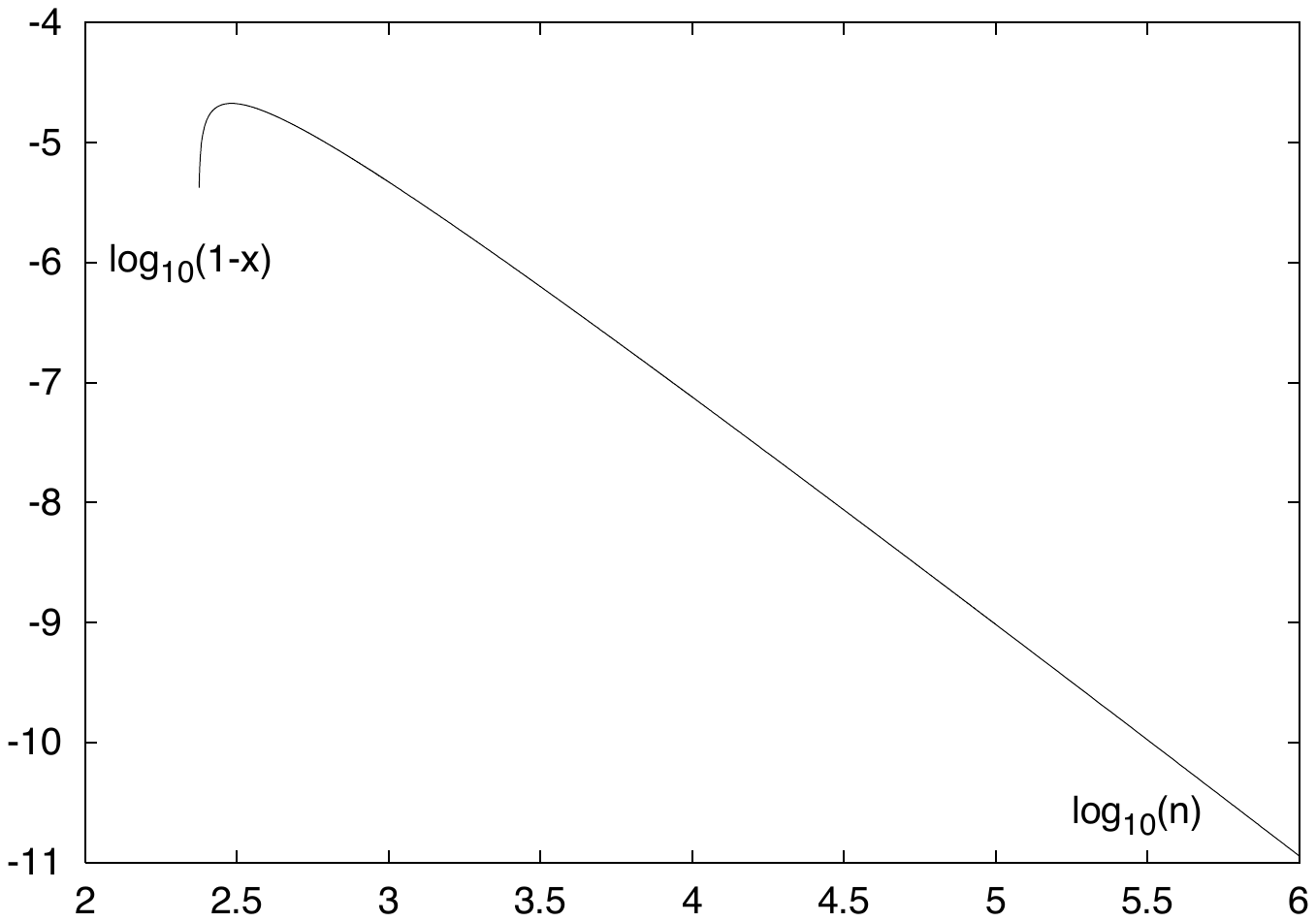}
\end{center}

\vspace*{2mm}

The compatibility of the system is equivalent to the proportionality of $(x^3-M, 1-M)$ to the first column $(\delta_n,k_n(x,\pi/n))$ of
$A_{n;2}$, that is

\vspace*{2mm}

$$M=\frac{x^3k_n(x,\fpn)-\delta_n}{k_n(x,\fpn)-\delta_n}.$$

\vspace*{2mm}

\noindent Finally, one must check the existence of positive masses. From system
$(*)$ it is easy to obtain numerically the solution as a one-dimension linear
manifold and to check that it intersects the positive octant. 

\vspace*{2mm}

\noindent Furthermore, for $n<473$ one has that solutions with $m_0=0$ are
admissible. This gives the really perverse solutions with minimal
$N$ known to present.

\vspace*{2mm}

\noindent Another interesting property is that for $n\le 1164$ the union
of the two regular polygons of a solution is a convex polygon with $2n$ vertices.

\vspace*{4mm}
\pagebreak

\begin{center} STUDY OF THE CASES OF $3$ POLYGONS \end{center}

\vspace*{2mm}

\noindent We shall consider two subcases

\vspace*{3mm}

\noindent {\bf 1) 2 equal radii: $\rho_2=\rho_3$}

\vspace*{3mm}

Assume $2n$ masses of values alternatively $m_2$ and $m_3$ are at the vertices
of a regular $2n$-gon centered at a mass $m_0.$ 

\vspace*{2mm}

\noindent Adding a regular $n$-gon also centered at $m_0$ and with all vertices
having mass $m_1$ we try to obtain a relative equilibrium giving rise to a
really perverse solution.

\vspace*{2mm}

\noindent It is not restrictive to assume $m_2\ge m_3$ and to take the radius
of the $2n$-gon equal to 1. Let $\rho$ be the radius of the $n$-gon. 

\vspace*{2mm}

\noindent By symmetry reasons, the polygon $P(m_1,n)$ must be homothetic either
to $P(m_2,n)$ or to $P(m_3,n).$ The two cases are studied separately. Our goal
is to prove the

\vspace*{1mm}

\begin{theo}[see \cite{B} Propositions 29 and 31] \label{t1} For any $n\ge 237$ there exists a really perverse
relative equilibrium consisting of a central mass $m_0$, a regular $n$-gon
with vertices of masses $m_1$ and a regular $2n$-gon with vertices of masses
alternatively equal to $m_2$ and $m_3$, the two polygons having the same
symmetry axes. All the masses can be chosen to be positive. For $n\le 236$ no really perverse solution of this type exists.
\end{theo}

\vspace*{1mm}

\noindent {\bf Proof.} We consider first the case $P(m_1,n)$ homothetic to
$P(m_2,n)$\footnote{The result is analogous in case $P(m_1,n)$ is homothetic to $P(m_3,n)$, the only difference being that now
$\rho>1$, that is: the $n$-gon is outside the $2n$-gon (see \cite{B})}. 
The relative equilibrium equations are
$$\begin{array}{rcrcrcr}
\rho^3-M&=&\delta_n\,nm_1&+&k_n(\rho^{-1},0)\,nm_2&+&k_n(\rho^{-1},\fpn)\,nm_3\\
\     1-M&=&k_n(\rho,0) \,nm_1&+&\delta_n \,nm_2&+&k_n(1,\fpn)\, nm_3\\
\     1-M&=&k_n(\rho,\fpn) \,nm_1&+&k_n(1,\fpn) \,nm_2&+&\delta_n \,nm_3
\end{array}$$

\noindent Let $\alpha_n(\rho)=\displaystyle{\frac{k_n(1,\fpn)-\delta_n}
{k_n(\rho,0)- k_n(\rho,\fpn)}}.$

\vspace*{3mm}

\noindent From the last two equations we obtain $ m_1=\alpha_n(\rho)(m_2-m_3).$ 
Hence $m_1=0$ if $m_2=m_3$. In particular, a relative equilibrium of the couple of a $2n$ gon and an $n$ gon, each one with with equal masses, cannot be really perverse while this is possible if the masses of the $2n$ gon are alternatively $m_2$ and $m_3\not=m_2$. 
\vspace*{3mm}

Before continuing with the proof we collect several useful informations.

\begin{lema}\label{lcolec} The following relations hold
\begin{itemize}
\item $k_n(1,\fpn)$ increases with $n.$ $\quad k_n(1,\fpn)<0$ iff $n\le 118$,
\item $k_n(1,\fpn)+\delta_n=2\delta_{2n}<0$ iff $n\le 236$,
\item $k_n(1,\fpn)-\delta_n=2(\delta_{2n}-\delta_n)>0$ for all $n\ge 2$,
\item $k_n(\rho,0)$ increases with respect to $\rho$ for all $n$. Furthermore
$k_n(0,0)=0,\; k_n(1^-,0)=+\infty,\; k_n(1^+,0)=-\infty,\; k_n(+\infty,0)=-1$,
\item $k_n(\rho,0)+k_n(\rho,\fpn)=2k_{2n}(\rho,0)> 0$ if $\rho\in[0,1)$, and
$< -2$ if $\rho>1$,
\item $k_n(\rho,0)-k_n(\rho,\fpn)\ge 0$ if $\rho\in[0,1)$ and
$k_n(\rho,0)-k_n(\rho,\fpn)\le 0$ if $\rho>1$,
\item $\alpha_n(\rho)$ decreases in $(0,1)\cup(1,+\infty).$ It can be made
continuous by defining $\alpha_n(1)=0.$ Furthermore $\alpha_n(0^+)=+\infty$
and $\alpha_n(+\infty)=-\infty.$
\end{itemize}   
\end{lema} 

\vspace*{2mm}

\noindent The proof of Lemma \ref{lcolec} follows from the integral
representation of $V_n(\rho,\theta)$ and from the expression of $\delta_n.$
Note that the fact that, while the limits of $k_n(\rho,0)$ when $\rho$ tends to 1 from the left or from the right are infinite, the value of $k_n(1,0)$ is well defined and finite. This is because the limit would correspond to adding the attraction
of the body at $\theta=0$ by a body coming from $x=1, l=n$ which we
excluded in the definition of $H_n(1,0)$.

.\smallskip

From the properties of $\alpha_n$ and $m_1>0,\; m_2\ge m_3$ it follows that
one should have $\rho<1$ (the $n$-gon is inside the $2n$-gon).  

\vspace*{-1mm}

\noindent Expressing $(1-\rho^3)M$ on the one hand by substracting the first
equilibrium equation to the product of the second by $\rho^3$, on the second
hand by replacing $m_1$ by $\alpha_n(\rho)(m_2-m_3)$ in the definition of $M$,
one obtains the equation

\vspace*{-3mm}

$$ \begin{array}{rcl}
0=f(\rho):&=&\fm (m_0+n\alpha_n(\rho)(m_2-m_3)+nm_2+nm_3)(1-\rho^3)\\ \
        & & + (\delta_n-\rho^3 k_n(\rho,0))\alpha_n(\rho)n(m_2-m_3)\\ \
        & & + (k_n(\rho^{-1},0)-\rho^3 \delta_n)nm_2+(k_n(\rho^{-1},\fpn)-\rho^3               k_n(1,\fpn))nm_3. \end{array} $$ 

\vspace*{1mm}

\noindent Looking for the dominant terms near the boundaries of $(0,1)$ and
using Lemma \ref{lcolec} we have $f(0^+)=+\infty,\quad f(1^-)=-\infty.$  

\vspace*{1mm}

\noindent Hence, for all $m_0,\,m_2\ge m_3,\, n$ there exist $\rho_0\in(0,1),\,
m_1>0$ giving rise to a relative equilibrium of the searched type.

\vspace*{2mm}

We now determine the conditions for such a relative equilibrium to be a really
perverse solution. The equations can be written as

\vspace*{-2mm}

$$ \left(\!\!\begin{array}{cc}
    k_n(\rho^{-1},0)\!+\!\delta_n\alpha_n(\rho) &
       k_n(\rho^{-1},\fpn)\!-\!\delta_n\alpha_n(\rho)\\ 
    \delta_n\!+\!k_n(\rho,0)\alpha_n(\rho) &  
       k_n(1,\fpn)\!-\!k_n(\rho,0)\alpha_n(\rho)
   \end{array}\!\!\right)\!
   \left(\!\!\begin{array}{c} nm_2\\ nm_3 \end{array}\!\!\right)\!=\!
   \left(\!\!\begin{array}{c}\rho^3\!-\!M\\  1\!-\!M \end{array}\!\!\right) $$

\noindent and, after adding the second column to the first one, the determinant
can be written

\vspace*{-2mm}

$$ \begin{array}{rcl}
d_n(\rho)&=&\left(k_n(\rho^{-1},0)+k_n(\rho^{-1},\fpn)\right)\left(
k_n(1,\fpn)-k_n(\rho,0)\alpha_n(\rho)\right)-\\ \
 & &\left(k_n(\rho^{-1},\fpn)-\delta_n\alpha_n(\rho)\right)\left(\delta_n+
k_n(1,\fpn)\right).\end{array}$$

\vspace*{2mm}

\noindent The dominant terms near the boundaries give

\vspace*{1mm}

$$ d_n(0^+)=\delta_n\left(\delta_n+k_n(1,\fpn)\right)\infty,\qquad
d_n(1^-)=-\delta_n\infty.$$

\vspace*{1mm}

\noindent Therefore, if $\delta_n+k_n(1,\fpn)>0$, which holds as soon a $n>236$ (recall that this function is increasing with $n$), 
there exists $\rho$ such that $d_n(\rho)=0.$
In fact, a somewhat long analysis
of the sign of $d_n(\rho)$ when $\delta_n+k_n(1,\fpn)<0$ (i.e., $n\le 236$),
making repeated use of Lemma \ref{lcolec} shows that it is positive for all
$\rho\in(0,1)$ (see \cite{B}).

\noindent The next figure shows the solution $\rho$ as a function of $n$. For
$n\in [237,472]$ the solution is close to zero. Going from 472 to 473 the
solution jumps to $\rho\approx 1.$ Then it remains close to 1 for $n\ge 473.$
For $n$ large it is better to use logarithmic scales.

\noindent It has been checked numerically that the solution whose existence
is proved, is unique.

\vspace*{3mm}

\noindent To study the compatibility of the equations defining a perverse solution we find the expression for $M$ obtained in writing that 
$(\rho^3-M, 1-M)$ is proportional to the sum of the columns of the matrix, that is

\vspace*{1mm}

$$M=\frac{\rho^3(\delta_n+k_n(1,\fpn))-(k_n(\rho^{-1},0)+k_n(\rho^{-1},\fpn))}
{\delta_n+k_n(1,\fpn)-(k_n(\rho^{-1},0)+k_n(\rho^{-1},\fpn))}.$$

\noindent As $\delta_n+k_n(1,\fpn)> 0,\; k_n(\rho^{-1},0)+k_n(\rho^{-1},\fpn)
<0$ (see Lemma \ref{lcolec}) the system is compatible.

\vspace*{4mm}

Finally we are interested in the positivity of the masses. After some 
manipulation of the equations (see \cite{B}) we obtain 

$$\begin{array}{cc} 0=m_0(1-\rho^3)+ K(\rho) \Delta m+ \\ \\
 n\left(k_n(\rho^{-1},0)+k_n(\rho^{-1},\fpn)
-\rho^3(\delta_n+k_n(1,\fpn))+2(1-\rho^3)\right)m_3,\end{array}$$

\vspace*{2mm}

\noindent where $K(\rho)$ is a function depending on $\rho$ and $\Delta m$
denotes $m_2-m_3.$ The term $m_0(1-\rho^3)$ is positive, $K(\rho) \Delta m$
can be chosen arbitrarily small by taking $\Delta m$ small enough. The
coefficient of $m_3$ is negative. Indeed $k_n(\rho^{-1},0)+k_n(\rho^{-1},\fpn)
+2\le 0$ and $\delta_n+k_n(1,\fpn)>0.$ Hence it is possible to have all the
masses positive.

\begin{center}

  \includegraphics[scale=0.65]{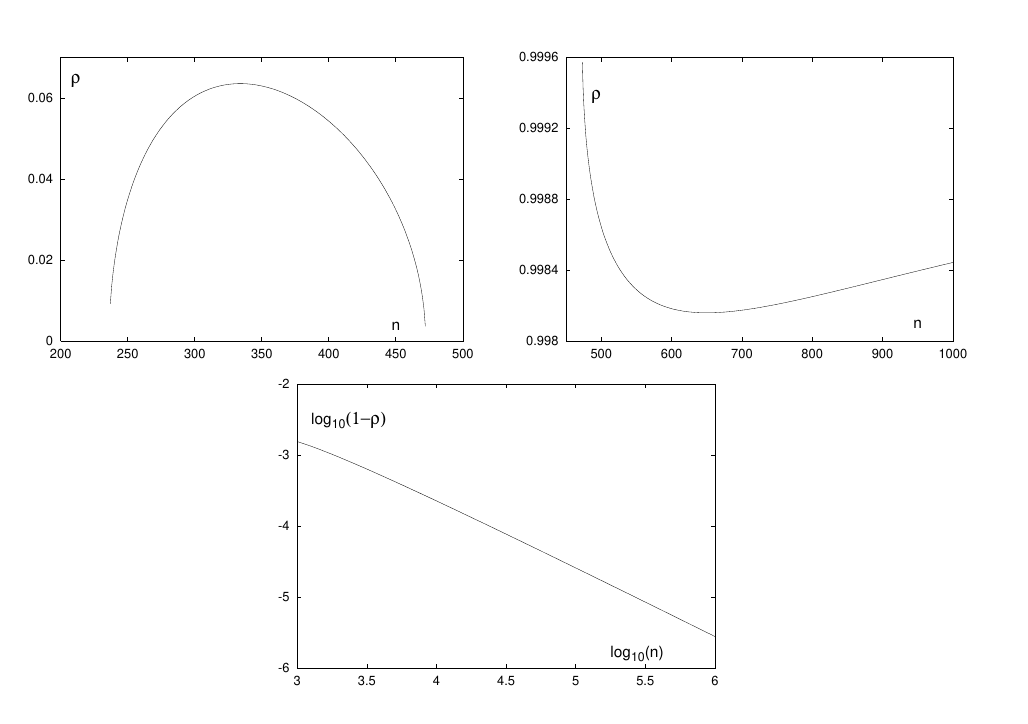}
\end{center}

\noindent {\bf 2) All radii different}

\vspace*{1mm}

\begin{prop} \label{p3p} For any integer $n$ between 456 and 874, there
exists a really perverse relative equilibrium solution of a central
mass surrounded by three homothetic regular $n$-gons. For any such $n$ there
is a closed curve of solutions, but not all points in these curves correspond
to positive masses. \end{prop}

\vspace*{1mm}

\noindent {\bf Proof.} Let $x_j=\rho_j/\rho_3$. One has to understand the rank 
of the matrix

$$ A_{n;3}(x_1,x_2)=\left(\begin{array}{ccc} 
\delta_n& k_n(x_2/x_1)&k_n(x_1^{-1})\\
k_n(x_1/x_2)& \delta_n& k_n(x_2^{-1})\\
k_n(x_1) & k_n(x_2) &\delta_n \end{array}\right),$$

\vspace*{1mm}

\noindent where, as before, $k_n(x)$ stands for $k_n(x,0).$

\vspace*{1mm}

\noindent We shall use the following notation:

$$ \hspace*{-4mm}\begin{array}{rcl}
 \det A_{n;3}(x_1,x_2)&\!=\!&\delta_n^3+p_n(x_1,x_2)\delta_n+q_n(x_1,x_2),\\ \
 p_n(x_1,x_2)&\!=\!& -k_n(x_1)k_n(\frac{1}{x_1})\!-\!k_n(x_2)
   k_n(\frac{1}{x_2})\!-\!k_n(\frac{x_1}{x_2})k_n(\frac{x_2}{x_1}),\\ \
 q_n(x_1,x_2)&\!=\!&k_n(x_1)k_n(\frac{1}{x_2})k_n(\frac{x_2}{x_1})+
      k_n(x_2)k_n(\frac{1}{x_1})k_n(\frac{x_1}{x_2}).	
   \end{array} $$

\vspace*{2mm}

\begin{lema}\label{lsim} $  \det A_{n;3}(x,x^\alpha)=\det A_{n;3}
(x,x^{1\!-\!\alpha})$, $\;p_n(x,x^\alpha)=p_n(x,x^{1\!-\!\alpha})$ and 

$\!\!\!\!\!q_n(x,x^\alpha)=q_n(x,x^{1-\alpha})$. \end{lema}

\vspace*{1mm}

\noindent {\bf Proof of the Lemma}. Replacing $\alpha$ by $1-\alpha$ transforms 
the matrix $A_{n;3}(x,x^\alpha)$ into its symmetrical one with respect to the second
diagonal.

\vspace*{4mm}

\noindent This lemma shows that, in the domain $0<x_1<x_2<1$, the functions
$\det A_{n;3},\, p_n$ and $q_n$  possess a symmetry which leaves fixed the curve
$x_2=\sqrt{x_1}$.

\begin{lema} \label{lboundary} The function $\det A_{n;3}$ has the same sign as
$\delta_n$ in some neighborhood of the triangular boundary of the region
$0<x_1<x_2<1$. Hence, it is negative in this neighborhood if
$\, n\le 472$ and positive if $n\ge 473$. \end{lema}

\vspace*{3mm}

\noindent The proof follows by an analysis of the dominant terms in different
pieces of the boundary.

\vspace*{3mm}

\begin{lema} \label{lroot} In the Newtonian case if $456\le n\le 472$, the
function $A_n(x,\sqrt{x})$ vanishes at least once on each of the intervals
$(0,0.098)$ and $(0.098,1).$ If $473\le n\le 874$, the function  
$A_n(x,\sqrt{x})$ vanishes at least once on each of the intervals
$(0,0.97)$ and $(0.97,1).$ \end{lema}

\vspace*{1mm}

\noindent {\bf Proof of Lemma \ref{lroot}} We must study the roots in the 
interval $(0,1)$ of the equation

$$\det A_{n;3}(x,\sqrt{x})=\delta_n^3+p_n(x,\sqrt{x})\delta_n+q_n(x,\sqrt{x})=0,
$$

\noindent where

$$\begin{array}{rcl}
p_n(x,\sqrt{x})&=&-k_n(x)k_n(1/x)-2k_n(\sqrt{x})k_n(1/\sqrt{x})\\ \
q_n(x,\sqrt{x})&=&k_n(x)k_n(1/\sqrt{x})^2+k_n(1/x)k_n(\sqrt{x})^2.
\end{array} $$

\vspace*{2mm}

\noindent From the analysis of $p_n(x,\sqrt{x})$ and $q_n(x,\sqrt{x})$ for
$x$ close to 0 and $x$ close to 1 it follows that if $n\le 472$ (resp. $n\ge 
473$) the determinant $A_{n;3}(x,\sqrt{x})$ is negative (resp. positive)
in a neighborhood of $x=0$ and $x=1$.

\vspace*{2mm}

\noindent Hence, to prove the Lemma, it is enough to check numerically
that
\begin{itemize}
\item for each integer $456\le n\le 472$, $A_{n;3}(0.098,\sqrt{0.098})$
is positive.
\item for each integer $473\le n\le 874$, $A_{n;3}(0.97,\sqrt{0.97})$ is 
negative.
\end{itemize}
This amounts to a finite number of checks.
\medskip

The question of the positivity of the masses is more involved. In fact, not
all the $(x_1,x_2)$ for which $\det A_{n;3}(x_1,x_2)=0$ have corresponding
masses such that all of them are positive.

\smallskip

\noindent The  figure below displays curves in $(x_1,x_2)$  for which
$\det A_{n;3}(x_1,x_2)=0$ for the admissible range of $n.$ The curve for
$n=456$ looks like a small ellipse around $P=(0.098,\sqrt{0.098})$, while
the one for $n=874$ looks like a small ellipse around $Q=(0.97,\sqrt{0.97})$.
The last curve of index 1 around $P$ has $n=472$, while the first curve
of index 1 around $Q$ has $n=473.$ 
\begin{center}
 \includegraphics[scale=0.45]{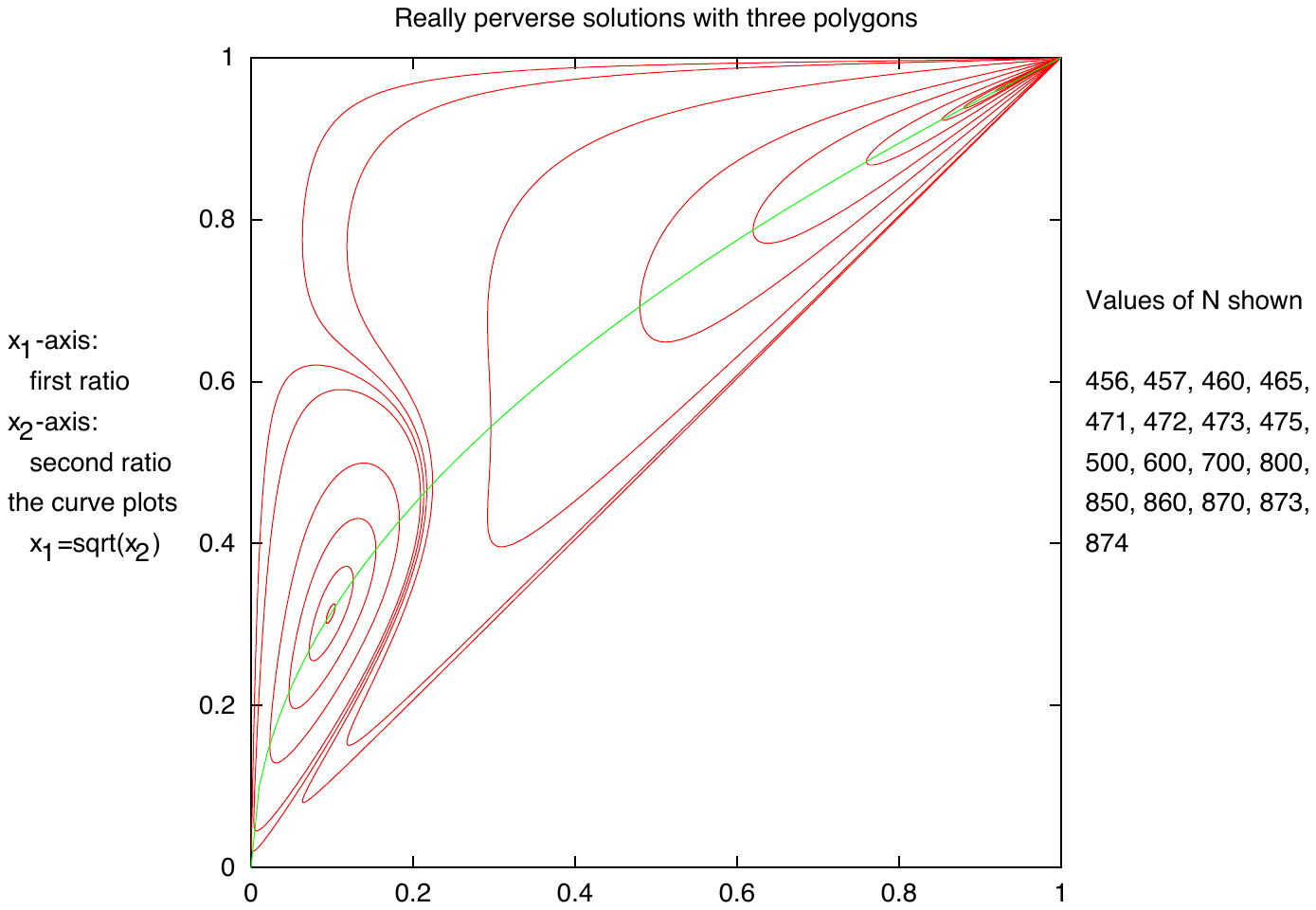}
\end{center}

\noindent The next figure magnifies the curves for
$n=870, 873$ and $874$. They are similar to very narrow ellipses.
\begin{center}
 \includegraphics[scale=0.3]{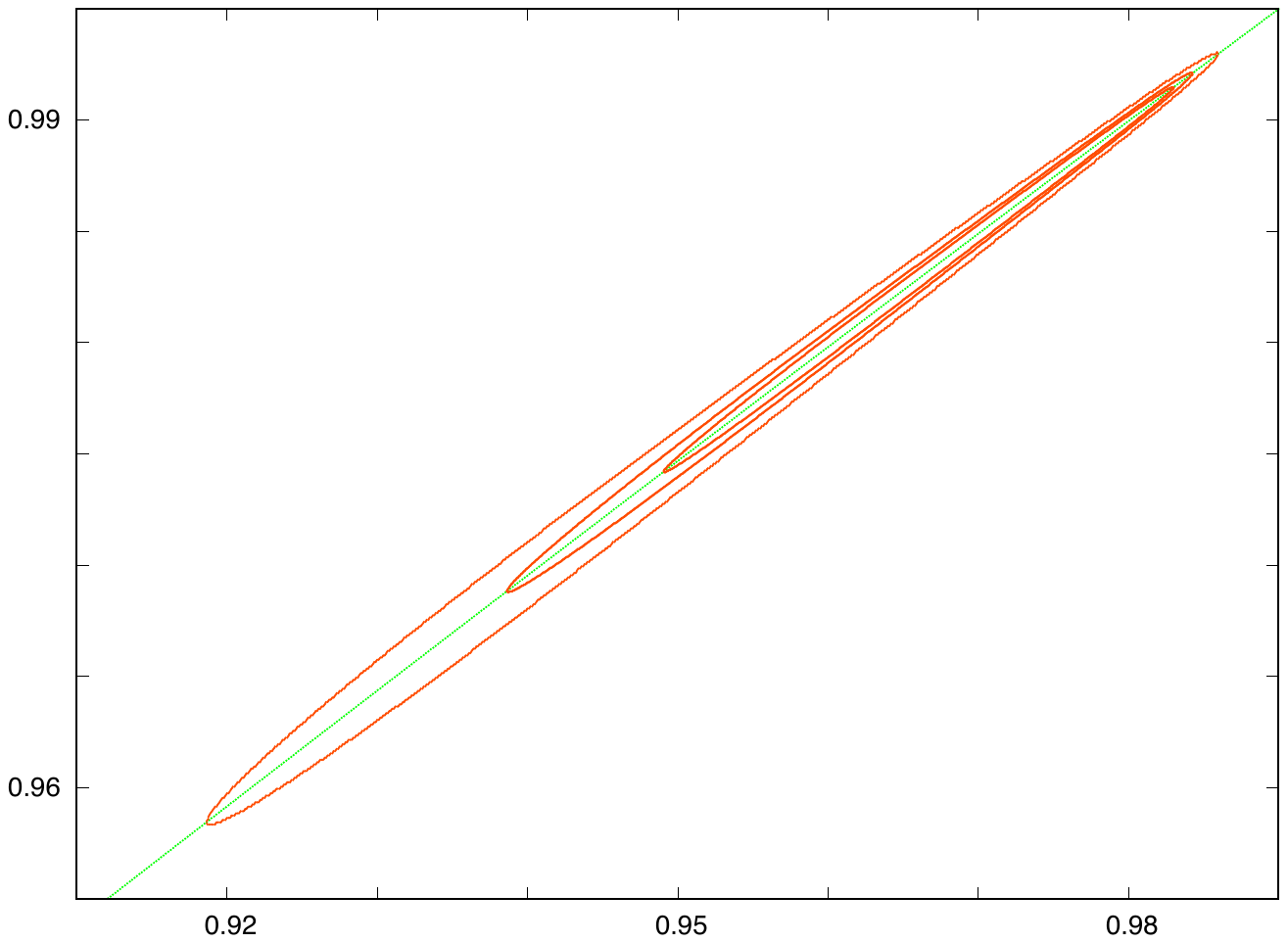}
\end{center}

\noindent In the third figure, beyond the curves for $n=472$ (on the left)
and $n=473$ (on the right), the plotted arcs correspond to the pieces of the
curves for $473\le n\le 609$ for which no solution with positive masses appears.
For the reamining values of $n\in [456,874]$ the full curve admits solutions
with positive masses.
\begin{center}
 \includegraphics[scale=0.55]{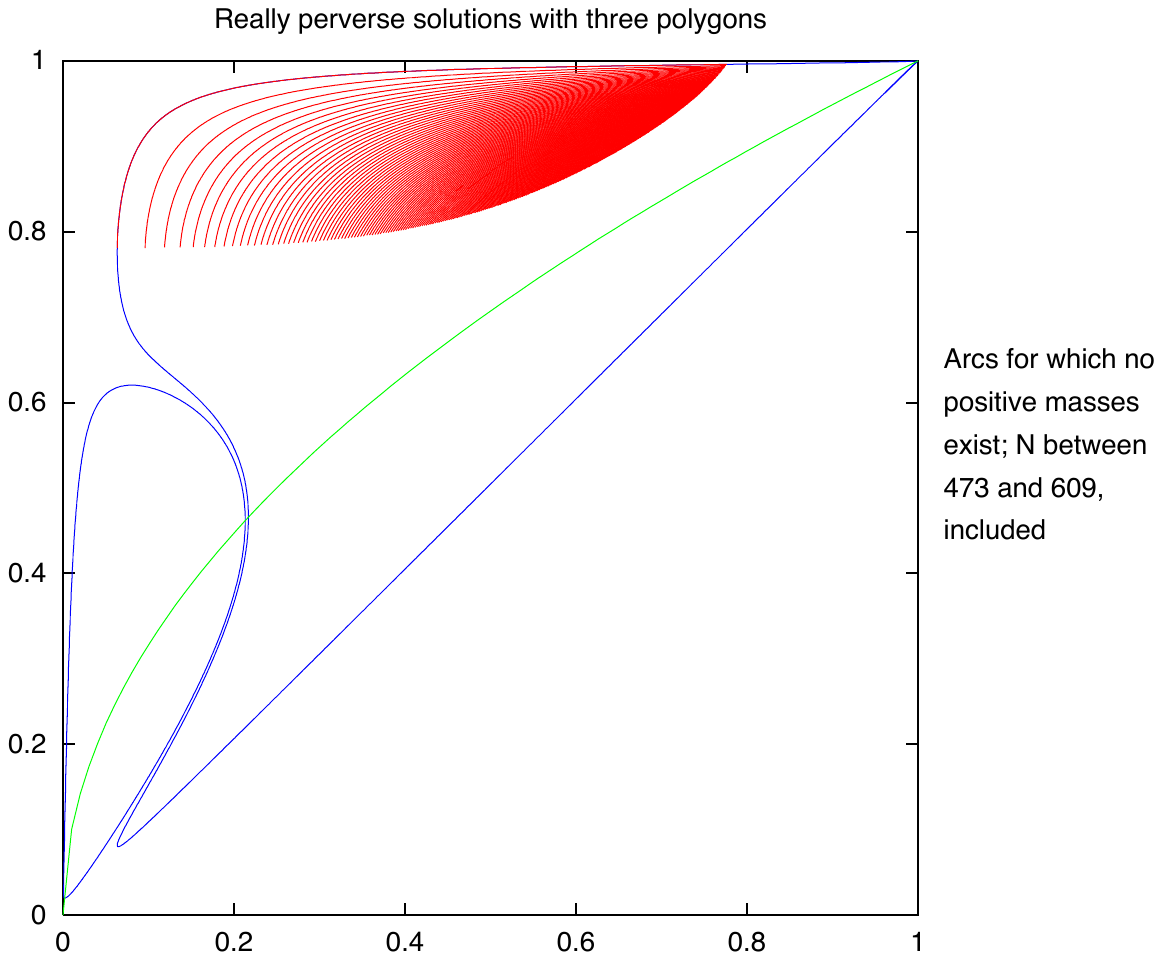}
\end{center}

\begin{center} THE CASE OF $k>3$ POLYGONS \end{center}

\vspace*{2mm}

The symmetry presented in Lemma \ref{lsim} persists in general. Let be
$x_j=\rho_j/\rho_k,\; j=1,\ldots,k-1$. We want to study the zeros of the
determinant of the matrix
$$ A_{n;k}(x_1,x_2,\ldots,x_{k\!-\!1})\!=\! \left(\!\!\begin{array}{ccccc}
\delta_n&k_n(\frac{x_2}{x_1})&\cdots&k_n(\frac{x_{k-1}}{x_1})&k_n(\frac{1}{x_1})
\\
 k_n(\frac{x_1}{x_2})&\delta_n&\cdots &\cdots& k_n(\frac{1}{x_2})\\
 \cdots&\cdots&\cdots&\cdots&\cdots\\
 k_n(\frac{x_1}{x_{k-1}})&\cdots&\cdots&\delta_n&k_n(\frac{1}{x_{k-1}})\\
 k_n(x_1)&k_n(x_2)&\cdots&k_n(x_{k-1})&  \delta_n \end{array}\!\!\right)\!. $$

\begin{lema} \label{lsimk} The following symmetry holds 
$$\det A_{n;k}(x,x^{\alpha_1},\ldots,x^{\alpha_{k-2}})=\det
A_{n;k}(x,x^{1-\alpha_{k-2}},\ldots,x^{1-\alpha_1}).$$ \end{lema}

\vspace*{2mm}

\!The proof is similar to the one of Lemma \ref{lsim}: replacing 
$\alpha_1,\ldots,\alpha_{k-2}$ by $1-\alpha_{k-2},\ldots,1\!-\!\alpha_1$,
respectively, transforms the initial matrix
$A_{n;k}(x,x^{\alpha_1},\ldots,x^{\alpha_{k-2}})$ into its symmetrical one with respect
to the second diagonal.

\vspace*{2mm}

We can introduce $\alpha_0=1,\;\alpha_{k-1}=0,$ that is 
$$x_1=x=x^{\alpha_0},\quad x_k=1=x^{\alpha_{k-1}}.$$

\vspace*{2mm}

\noindent {\bf Fixed points of the symmetry}
\begin{itemize}
\item [1)] $k$ odd:
$$ \alpha_1+\alpha_{k-2}=1,\ldots ,\alpha_{\frac{k-3}{2}}+\alpha_{\frac{k+1}{2}}
=1,\; \alpha_{\frac{k-1}{2}}=\frac{1}{2}.$$
\item [2)] $k$ even:
$$\alpha_1+\alpha_{k-2}=1,\ldots,\alpha_{\frac{k-2}{2}}+\alpha_{\frac{k}{2}}=1.$$
\end{itemize}

\vspace*{2mm}

The existence of really perverse relative equilibria with more than 3 polygons
is a direct consequence of the existence for 3 polygons. Consider, for instance,
ratios of radii of the form $(\eps^{k-3},\ldots,\eps^2,\eps,x_1,x_2)$, where
$(x_1,x_2)$ are ratios for the case $k=3.$ 

\vspace*{2mm}

But there are other possibilities. If one has a solution with $k=3$ and
some $(\rho_1,\rho_2,\rho_3=1)$, then there are related solutions with 
$k=4$ of the form $(\eps,\rho_1,\rho_2,\rho_3)$ if $\eps$ is small and small
changes are done to $(\rho_1,\rho_2,\rho_3)$. Another possibility is
$(\eps\rho_1,\eps\rho_2,\eps\rho_3,1)$.

\vspace*{2mm}

\noindent {\bf All these kinds of solutions have been found numerically}.

\vspace*{2mm}

\noindent With the same idea, for $k=5$ and if $\rho_1,\rho_2,\rho_3$ are
suitable values for $k=3$, one can take

$(\eps^2, \eps,\rho_1,\rho_2,\rho_3)$, but also

$(\eps^2, \eps\rho_1,\eps\rho_2,\eps\rho_3,1)$ or

$(\eps^2\rho_1,\eps^2\rho_2,\eps^2\rho_3,\eps,1)$.

\vspace*{2mm}

\noindent For $k=6$ there are different places where to put powers of $\eps$:

$(\eps^3, \eps^2, \eps, \rho_1,\rho_2,\rho_3)$,

$(\eps^3, \eps^2, \eps\rho_1,\eps\rho_2,\eps\rho_3,1)$,

$(\eps^3, \eps^2\rho_1,\eps^2\rho_2,\eps^2\rho_3,\eps,1)$,

$(\eps^3\rho_1,\eps^3\rho_2,\eps^3\rho_3,\eps^2,\eps,1)$,

\noindent and it is also possible, if beyond $\rho_1,\rho_2,\rho_3=1$ another
solution $\mu_1,\mu_2,$ $\mu_3=1$ is available for $k=3$, to form solutions like

$(\eps\rho_1,\eps\rho_2,\eps\rho_3,\mu_1,\mu_2,\mu_3)$.

\vspace*{4mm}

The range of values of $n$ for which really perverse solutions with $k$
homothetic $n$-gons exist depends also on $k$. For $k=5$ the range is
$451\le n\le 1265.$ It is worth to remark that for these extreme values
one finds surfaces like ellipsoides in $\Rset^4$ centered around points of
the form $\rho_1=\eps^4, \rho_2=\eps^3, \rho_3=\eps^2, \rho_4=\eps, \rho_5=1$.
For $n=451$ one has $\eps\approx 0.314$ and for $n=1265$ the value is
$\eps\approx 0.9873.$

Using $k$ odd and values $\rho_j=\eps^{k-j},\, j=1,\ldots,k$ the extremal
values of $n$ (minimal and maximal) seem to correspond (at least locally) 
to the extrema of $n$ with $k$ polygons. 

\vspace*{2mm}

Going up to $k=51$ the minimal $n$ is 442 and no further variation seems to
be expected. The corresponding value of $\eps$ tends to $0.316\ldots .$

\vspace*{2mm}

On the other hand, the maximal $n$ seems to increase without bound with $k$.
For $k=51$ it reaches the value $n_{max}(51)=9911$, and the derivative
$dn_{max}(k)/dk$ seems to stabilize near 373. The corresponding value of $\eps$
approaches 1. For $k=51$ it is $\eps_{max}=0.998879$.

\vspace*{4mm}

\begin{center} OTHER HOMOGENEOUS POTENTIALS \end{center}

\vspace*{4mm}

Instead of the Newtonian potential $r^{-1}$ one can consider potentials of
the form $r^{-2\beta}.$ Then all the functions in the analysis of the problem
depend also on $\beta.$ In particular $\delta_n=\delta_n(\beta)$ and the
transitions occur at different values of $n.$

\vspace*{2mm}

Concerning the behaviour of $\delta_n(\beta)$ with respect to $n$ for $n$
large, it is clear that $\delta_n(1/2)=O(\log(n))$ and $\delta_n(\beta)=
O(n^{2\beta-1})$ if $\beta>1/2.$ But for $\beta<1/2$ the values of
$\delta_n(\beta)$ are bounded when $n\to +\infty$.

\vspace*{2mm}

Given a value of $n>6$ there exist $\beta(n)$ such that $\delta_n(\beta(n))=0.$
The function $\beta(n)$ is decreasing with respect to $n.$ When $n\to +\infty$
$\beta(n)$ tends to a limit value $\beta_c=0.395\ldots\,.$ $\beta(n)$ does not
change sign if $\beta\le \beta_c.$

\begin{center}
 \includegraphics[scale=0.35]{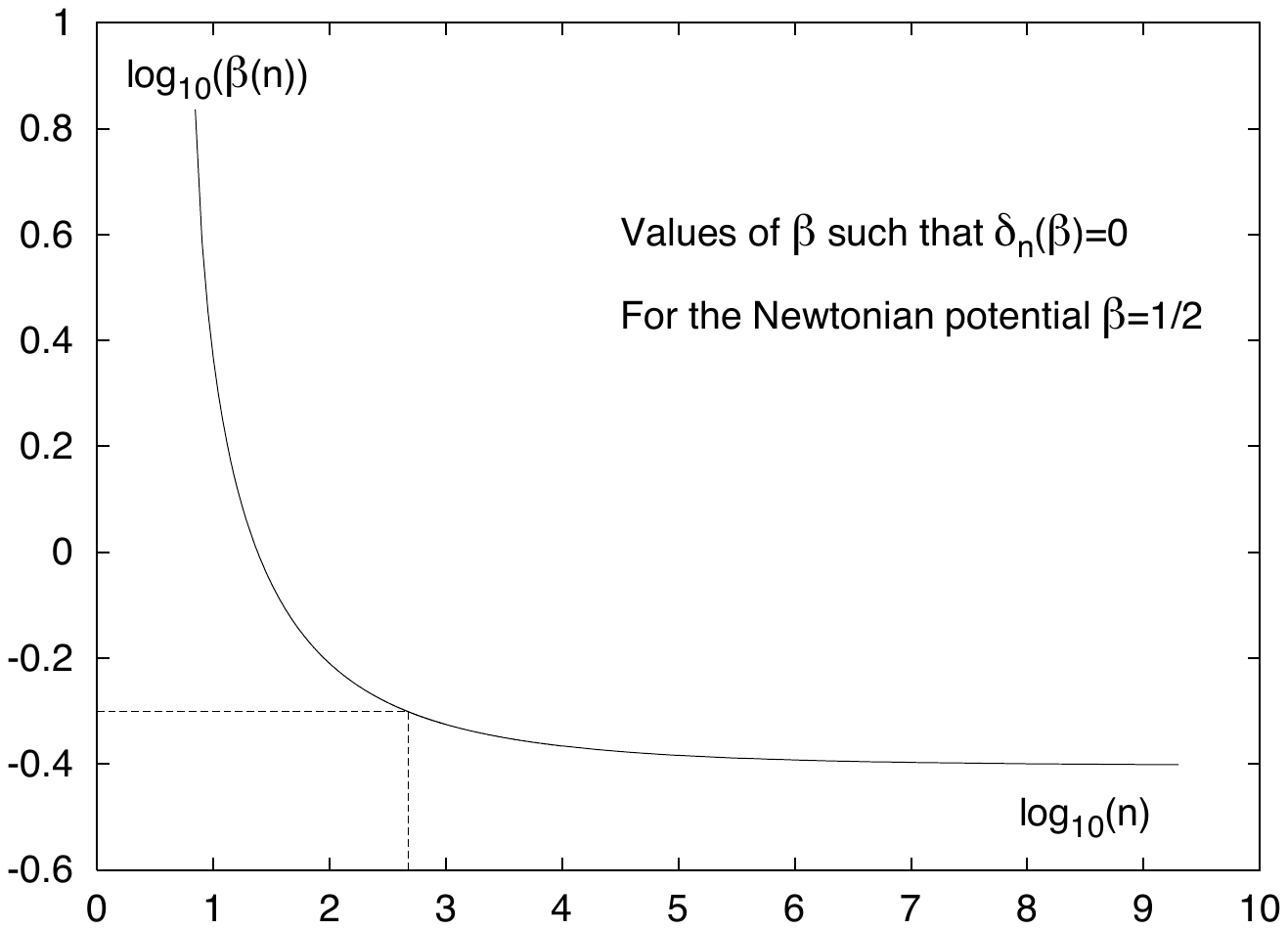}
\end{center}

On the other hand $\beta(7)= 6.85879\ldots$. For $\beta>\beta(7)$ the $\delta_n$
changes sign between $n=6$ and $n=7$. For these potentials the case of 2
shifted polygons could be realizable with $N=8$ if $m_0=0$ is admissible.

\noindent {\bf Additional note:} {\it These results were presented by C. Sim\'o at
the Workshop on Differential Equations dedicated to the memory of Vladimir F.
Lazutkin, in St. Petersburg, August 2002. Questions posed by Prof. M. Shibayama
have suggested the authors to make these results public, after revising the text
and updating references and affiliations. The authors are very grateful to both referees who pointed an error, now corrected, 
in the formula for $H_n$.}
}
\end{document}